\documentclass[11pt]{article}
\usepackage{fullpage} 
\usepackage{blkarray}
\usepackage{ucs}
\usepackage{amssymb}
\usepackage{amsthm}
\usepackage{amsmath}
\usepackage{latexsym}
\usepackage[utf8]{inputenc}
\usepackage[english]{babel}
\usepackage{graphicx}
\usepackage{wrapfig}
\usepackage{caption}
\usepackage[margin=0.81in]{geometry}
\usepackage{subcaption}
\usepackage{txfonts}
\usepackage{lmodern}
\usepackage{mathrsfs}
\usepackage{algorithm}
\usepackage{dsfont}
\usepackage{enumerate}
\usepackage[hidelinks]{hyperref}
\usepackage{multicol}
\usepackage{csquotes}
\usepackage{stmaryrd}
\usepackage{tikz}
\usepackage{comment}
\usepackage{enumitem}
\usepackage{cite}
\newtheorem{theo}{Theorem}

\newtheorem{lemma}[theo]{Lemma}
\newtheorem{ques}[theo]{Question}

\numberwithin{theo}{section}

\tikzset{
vtx/.style={inner sep=1.1pt, outer sep=0pt, circle, fill,draw}
}

\usepackage{xparse,xpatch,xcolor,tikz}
\usepackage{pgfmath}
\usepackage{xifthen}
\usetikzlibrary{decorations.text, arrows.meta,calc,shadows.blur,shadings}
\usetikzlibrary{shapes.geometric, arrows, positioning}

\author{%
Felix Christian Clemen \footnote {Karlsruhe Institute of Technology, 76133 Karlsruhe, Germany, E-mail: \texttt{felix.clemen@kit.edu}.}
 \and Adam Zsolt Wagner \footnote {Worcester Polytechnic Institute, Worcester, Massachusetts 01609, USA, E-mail: \texttt{zadam@wpi.edu}. }
}

\begin{document}

\title{A note on balanced edge-colorings avoiding rainbow cliques of size four}

\maketitle

\abstract{A balanced edge-coloring of the complete graph is an edge-coloring such that every vertex is incident to each color the same number of times. In this short note, we present a construction of a balanced edge-coloring with six colors of the  complete graph on $n=13^k$ vertices, for every positive integer $k$, with no rainbow $K_4$. This solves a problem by Erd\H{o}s and Tuza.}

\section{Introduction and main result}
Let $H$ be a graph. An edge-coloring of the complete graph $K_n$ contains a \emph{rainbow copy} of $H$ if it contains a copy of $H$ such that all its edges are assigned  different colors. Various conditions on edge-colorings forcing the existence of rainbow copies have been studied~\cite{AJMP,APS,AJT,ErdTu,KMSV,MBNL,RT,SS}. Here, we look at a condition where the colors are well distributed. A \emph{balanced} edge-coloring of the complete graph is an edge-coloring such that every vertex is incident to each color the same number of times.

\begin{ques}[Erd\H{o}s and Tuza, Problem 1 in~\cite{ErdTu}]
\label{problem1}
Is the following true for every graph $H$? For $n$ sufficiently large, every balanced edge-coloring  with $|E(H)|$ colors of the complete graph $K_n$ contains a rainbow copy of $H$.
\end{ques}

Recently, this question was answered by Axenovich and Clemen~\cite{AxCl} who showed that it is not true for most cliques. In particular, they showed that it is false for all cliques $H=K_q$ with odd number of edges and at least six vertices. They~\cite{AxCl} also conjectured Question~\ref{problem1} to be false for every clique with at least four vertices. Here, we answer Question~\ref{problem1} in the negative for $H=K_4$, a subcase which has been given particular attention by Erd\H{o}s and Tuza.

Indeed, Erd\H{o}s and Tuza~\cite{ErdTu} suggested that the simplest counterexamples to Question~\ref{problem1} might be $K_4$, $C_6$ and $2K_3$ and commented ``we could not prove or disprove that every $t$-regular $6$-coloring of $K_{6t+1}$, contains them as rainbow subgraphs (here $t$ has to be even)". In Erd\H{o}s' list ``Some of my favourite problems on cycles and colourings"~\cite{erdos1996some}, he further remarked that one of the most interesting unsolved problems in this area is Question~\ref{problem1} for $H=C_6$ and $H=K_4$.  

\begin{theo}
\label{K4rainbow}
For every $k\geq 1$ there exists a balanced edge-coloring of $K_{13^k}$ with $6$ colors and no rainbow $K_4$.
\end{theo}
A key idea in the constructions by Axenovich and Clemen~\cite{AxCl} is that iterating a balanced coloring with no rainbow copy of some clique $K_q$ maintains those properties: 
\begin{lemma}[Axenovich, Clemen, Lemma 2.2 in~\cite{AxCl}]
\label{iteralexico}
If there exists a balanced edge-coloring of $K_n$ with $\ell$ colors and no rainbow $K_q$, then for every $k\geq 1$ there exists a balanced edge-coloring of $K_{n^k}$ with $\ell$ colors and no rainbow $K_q$. 
\end{lemma}
Therefore, it suffices to prove Theorem~\ref{K4rainbow} for $k=1$. The base constructions used in \cite{AxCl} are the standard examples of 1- and 2-factorizations, i.e.~edge-colorings such that every color class is a 1-regular, respectively 2-regular, spanning subgraph. Note that those colorings do not work as constructions for Theorem~\ref{K4rainbow}. 

\begin{figure}[h]
    \centering
    \begin{subfigure}{0.45\textwidth}
    \centering
\begin{equation*}
\begin{pmatrix}

0 & 2 &5&4&1&3&3&6&4&2&6&5&1 \\
2&0&3&6&5&6&4&1&3&1&4&5&2 \\
5&3&0&5&4&2&6&3&1&6&2&1&4 \\
4&6&5&0&2&4&5&2&1&3&3&1&6 \\
1&5&4&2&0&3&1&6&2&5&4&6&3 \\
3&6&2&4&3&0&1&4&5&6&5&2&1 \\
3&4&6&5&1&1&0&2&5&4&2&6&3 \\
6&1&3&2&6&4&2&0&3&5&1&4&5 \\
4&3&1&1&2&5&5&3&0&4&6&2&6 \\
2&1&6&3&5&6&4&5&4&0&1&3&2 \\
6&4&2&3&4&5&2&1&6&1&0&3&5 \\
5&5&1&1&6&2&6&4&2&3&3&0&4 \\
1&2&4&6&3&1&3&5&6&2&5&4&0 
\end{pmatrix}
\end{equation*}
    \caption{Adjacency matrix of the coloring}
	\label{subfig:1}
    \end{subfigure}
    \begin{subfigure}{0.5\textwidth}
    \centering
\begin{tikzpicture}[scale=1.2]
\draw
\foreach \i in {1,...,13}{
(90+27.692*\i:2.5) coordinate(\i) node[vtx]{}
};
\draw[color=red] (1)  to (5);
\draw[color=red] (1)  to (13);
\draw[color=red] (2)  to (8);
\draw[color=red] (2)  to (10);
\draw[color=red] (3)  to (9);
\draw[color=red] (3)  to (12);
\draw[color=red] (4)  to (9);
\draw[color=red] (4)  to (12);
\draw[color=red] (5)  to (7);
\draw[color=red] (6)  to (7);
\draw[color=red] (6)  to (13);
\draw[color=red] (8)  to (11);
\draw[color=red] (10)  to (11);

\draw[color=blue] (1)  to (2);
\draw[color=blue] (1)  to (10);
\draw[color=blue] (2)  to (13);
\draw[color=blue] (3)  to (6);
\draw[color=blue] (3)  to (11);
\draw[color=blue] (4)  to (5);
\draw[color=blue] (4)  to (8);
\draw[color=blue] (5)  to (9);
\draw[color=blue] (6)  to (12);
\draw[color=blue] (7)  to (8);
\draw[color=blue] (7)  to (11);
\draw[color=blue] (9)  to (12);
\draw[color=blue] (10)  to (13);

\draw[color=green] (1)  to (6);
\draw[color=green] (1)  to (7);
\draw[color=green] (2)  to (3);
\draw[color=green] (2)  to (9);
\draw[color=green] (3)  to (8);
\draw[color=green] (4)  to (10);
\draw[color=green] (4)  to (11);
\draw[color=green] (5)  to (6);
\draw[color=green] (5)  to (13);
\draw[color=green] (7)  to (13);
\draw[color=green] (8)  to (9);
\draw[color=green] (10)  to (12);
\draw[color=green] (11)  to (12);

\draw[color=brown] (1)  to (4);
\draw[color=brown] (1)  to (9);
\draw[color=brown] (2)  to (7);
\draw[color=brown] (2)  to (11);
\draw[color=brown] (3)  to (5);
\draw[color=brown] (3)  to (13);
\draw[color=brown] (4)  to (6);
\draw[color=brown] (5)  to (11);
\draw[color=brown] (6)  to (8);
\draw[color=brown] (7)  to (10);
\draw[color=brown] (8)  to (12);
\draw[color=brown] (9)  to (10);
\draw[color=brown] (12)  to (13);

\draw[color=gray] (1)  to (3);
\draw[color=gray] (1)  to (12);
\draw[color=gray] (2)  to (5);
\draw[color=gray] (2)  to (12);
\draw[color=gray] (3)  to (4);
\draw[color=gray] (4)  to (7);
\draw[color=gray] (5)  to (10);
\draw[color=gray] (6)  to (9);
\draw[color=gray] (6)  to (11);
\draw[color=gray] (7)  to (9);
\draw[color=gray] (8)  to (10);
\draw[color=gray] (8)  to (13);
\draw[color=gray] (11)  to (13);

\draw[color=yellow] (1)  to (8);
\draw[color=yellow] (1)  to (11);
\draw[color=yellow] (2)  to (4);
\draw[color=yellow] (2)  to (6);
\draw[color=yellow] (3)  to (7);
\draw[color=yellow] (3)  to (10);
\draw[color=yellow] (4)  to (13);
\draw[color=yellow] (5)  to (8);
\draw[color=yellow] (5)  to (12);
\draw[color=yellow] (6)  to (10);
\draw[color=yellow] (7)  to (12);
\draw[color=yellow] (9)  to (11);
\draw[color=yellow] (9)  to (13);

\path (0,1.2) -- (0,-3);
\end{tikzpicture}
    \caption{A drawing of the coloring}
	\label{subfig:2}
    \end{subfigure}
    \caption{The edge-coloring of $K_{13}$ with no rainbow $K_4$ }
    \label{fig:3}
\end{figure}
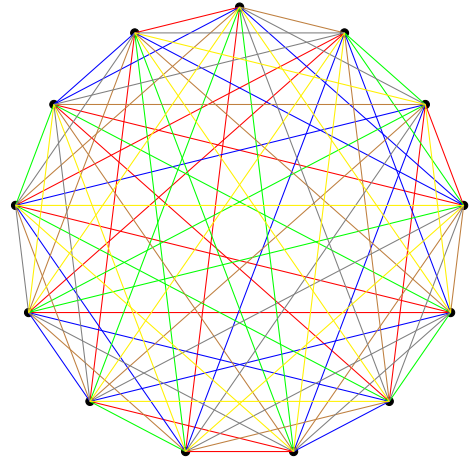

 In Figure~\ref{fig:3}, we present an edge-coloring of $K_{13}$ with six colors, found by a computer search, such that every vertex is incident to every color exactly twice and there is no rainbow $K_4$. While it can be checked quickly that this coloring is indeed balanced, it takes more effort to see that it does not contain a rainbow $K_4$, simply because there are $\binom{13}{4}=715$ copies of $K_4$ in $K_{13}$. 
 
 We remark that our construction differs from the constructions in \cite{AxCl} in the sense that it seemingly does not follow a visible pattern.

\vspace{-0.1cm}

\vspace{-0.5cm}
\section*{Acknowledgments}
The first author thanks Maria Axenovich for introducing him to Question~\ref{problem1} and Bernard Lidick\'{y}  for discussions on the topic.

\bibliographystyle{abbrv}
\bibliography{rainbowK4}

\begin{thebibliography}{10}

\bibitem{AJMP}
N.~Alon, T.~Jiang, Z.~Miller, and D.~Pritikin.
\newblock Properly colored subgraphs and rainbow subgraphs in edge-colorings
  with local constraints.
\newblock {\em Random Structures \& Algorithms}, 23(4):409--433, 2003.

\bibitem{APS}
N.~Alon, A.~Pokrovskiy, and B.~Sudakov.
\newblock Random subgraphs of properly edge-coloured complete graphs and long
  rainbow cycles.
\newblock {\em Israel J. Math.}, 222(1):317--331, 2017.

\bibitem{AxCl}
M.~Axenovich and F.~C. Clemen.
\newblock Rainbow subgraphs in edge-colored complete graphs - answering two
  questions by {E}rd{\H{o}}s and {T}uza.
\newblock {\em arXiv:2209.13867}, 2022.

\bibitem{AJT}
M.~Axenovich, T.~Jiang, and Z.~Tuza.
\newblock Local anti-{R}amsey numbers of graphs.
\newblock {\em Combin. Probab. Comput.}, 12(5-6):495--511, 2003.
\newblock Special issue on Ramsey theory.

\bibitem{erdos1996some}
P.~Erd\H{o}s.
\newblock Some of my favourite problems on cycles and colourings.
\newblock {\em Tatra Mt. Math. Publ}, 9:7--9, 1996.

\bibitem{ErdTu}
P.~Erd\H{o}s and Z.~Tuza.
\newblock Rainbow subgraphs in edge-colorings of complete graphs.
\newblock {\em Ann. Discrete Math.}, 55:81--88, 1993.

\bibitem{KMSV}
P.~Keevash, D.~Mubayi, B.~Sudakov, and J.~Verstra\"{e}te.
\newblock Rainbow {T}ur\'{a}n problems.
\newblock {\em Combin. Probab. Comput.}, 16(1):109--126, 2007.

\bibitem{MBNL}
J.~J. Montellano-Ballesteros and V.~Neumann-Lara.
\newblock An anti-{R}amsey theorem.
\newblock {\em Combinatorica}, 22(3):445--449, 2002.

\bibitem{RT}
V.~R\"{o}dl and Z.~Tuza.
\newblock Rainbow subgraphs in properly edge-colored graphs.
\newblock {\em Random Structures \& Algorithms}, 3(2):175--182, 1992.

\bibitem{SS}
M.~Simonovits and V.~T. S\'{o}s.
\newblock On restricted colourings of {$K_n$}.
\newblock {\em Combinatorica}, 4(1):101--110, 1984.

\end{thebibliography}

\end{document}